\documentclass[runningheads]{llncs}

\usepackage{algorithm}
\usepackage{algpseudocode}
\usepackage{amsmath,amssymb,amsfonts}
\usepackage{array}
\usepackage{multirow}
\usepackage{microtype}
\usepackage{graphicx}
\usepackage{tabularx}
\usepackage[utf8]{inputenc}
\usepackage[english]{babel}

\DeclareMathOperator*{\argmin}{arg\,min}

\newcolumntype{Y}{>{\centering\arraybackslash}X}

\begin{document}
    \title{Lossy Image Compression with Stochastic Quantization}

    \author{Anton Kozyriev\inst{1} \and Vladimir Norkin\inst{1,2}}
    \authorrunning{Anton Kozyriev and Vladimir Norkin}

    \institute{National Technical University of Ukraine ''Igor Sikorsky Kyiv Polytechnic Institute'', Kyiv, 03056, Ukraine \email{a.kozyriev@kpi.ua} \and
    V.M.Glushkov Institute of Cybernetics, National Academy of Sciences of Ukraine, Kyiv, 03178, Ukraine \email{v.norkin@kpi.ua}}

    \maketitle

    \begin{abstract}

    Lossy image compression algorithms play a crucial role in various domains, including graphics, and image processing. As image information density increases, so do the resources required for processing and transmission. One of the most prominent approaches to address this challenge is color quantization, proposed by Orchard et al. (1991). This technique optimally maps each pixel of an image to a color from a limited palette, maintaining image resolution while significantly reducing information content. Color quantization can be interpreted as a clustering problem (Krishna et al. (1997), Wan (2019)), where image pixels are represented in a three-dimensional space, with each axis corresponding to the intensity of an RGB channel. However, scaling of traditional algorithms like K-Means can be challenging for large data, such as modern images with millions of colors. This paper reframes color quantization as a three-dimensional stochastic transportation problem between the set of image pixels and an optimal color palette, where the number of colors is a predefined hyperparameter. We employ Stochastic Quantization (SQ) with a seeding technique proposed by Arthur et al. (2007) to enhance the scalability of color quantization. This method introduces a probabilistic element to the quantization process, potentially improving efficiency and adaptability to diverse image characteristics. To demonstrate the efficiency of our approach, we present experimental results using images from the ImageNet dataset. These experiments illustrate the performance of our Stochastic Quantization method in terms of compression quality, computational efficiency, and scalability compared to traditional color quantization techniques.

    \end{abstract}

    \section{Introduction}

    Digital image representation in modern video displays relies on additive color mixing, where the intensity of three primary colors (Red, Green, and Blue) is modulated at each pixel \cite{Orchard_Bouman_1991,sharma2017digital}. In this system, each pixel is encoded as a triplet of unsigned integers, with each element corresponding to the intensity of a primary color. For instance, an 8-bit display can represent $2^8 = 256$ values per RGB channel, yielding a total of $2^{8 \times 3} = 16,777,216$ possible colors. To efficiently store images produced by digital photography, compression algorithms have been developed to balance data storage requirements and image quality. Many of these algorithms fall into the category of lossy compression, where the compressed image experiences a reduction in quality that cannot be fully recovered to its original state.

    Color quantization, introduced by Orchard et al. \cite{Orchard_Bouman_1991}, is one approach to lossy image compression that reduces image data without compromising the original resolution. This method comprises two primary stages: selecting an optimal color palette of dominant colors present in the image, and mapping each pixel to the nearest color within this optimal palette. This approach conceptualizes each pixel as a point in three-dimensional space, with each axis representing the intensity of one of the three primary colors. Consequently, pixels with similar color shades cluster more closely in this space. The original paper explores various techniques to address the color quantization problem, including total squared error (TSE) minimization, binary tree palette design, and the Linde-Buzo-Gray (LBG) algorithm, among others. Krishna et al. \cite{Krishna_1997} proposed a solution to the color quantization problem utilizing the K-Means algorithm \cite{Lloyd_1982} for optimal color palette selection, where cluster centers represent optimal colors. This approach was further investigated in subsequent works \cite{Celebi_2011,Wan_2019}.

    However, a recent study \cite{Norkin_2024} highlighted the limitations of traditional clustering methods in solving quantization problems, particularly their poor scalability for large datasets. With the rapid increase in image resolution in modern video displays, this scalability bottleneck may significantly impact the performance of these color quantization methods.

    In this paper, we present a novel interpretation of color quantization within the domain of stochastic programming, specifically as a non-convex stochastic transportation problem \cite{Kuzmenko_Uryasev_2019}. We employ the Stochastic Quantization algorithm \cite{Lakshmanan_Pichler_2023,Norkin_2024} to determine the optimal color palette by minimizing the distance to the set of colors in the original image. We provide experimental results of this approach using test images from the ImageNet dataset \cite{Russakovsky_2015}, implemented using a Python implementation of the Stochastic Quantization algorithm.

    \section{The Problem Setting}

    The transportation problem (see, e.g., \cite{Kuzmenko_Uryasev_2019,Lakshmanan_Pichler_2023,Norkin_2024}) is utilized to approximate one discrete probability distribution $\{\xi_i\} \in \Xi \subset \mathbb{R}^{n}$ with another discrete distribution containing fewer elements $\{y_k\} \in Y \subset \mathbb{R}^{n}$. The optimal positioning of each element $\{y_1, \ldots, y_K\}$ is determined by minimizing the Wasserstein (Kantorovich–Rubinstein) distance:

    \begin{equation}
        \label{transport:eq}
            \min_{y = \{ y_1, \ldots, y_K \} \in Y^K \subset \mathbb{R}^{nK}} \min_{q = \{ q_1, \ldots, q_K \} \in \mathbb{R}^K_{+}} \min_{x = \{ x_{ij} \geq 0 \}} \sum_{i=1}^I \sum_{k=1}^K d(\xi_i, y_k)^r x_{ik}
    \end{equation}

    subject to the constraints:

    \begin{equation}
        \label{transport-constraints:eq}
            \sum_{k=1}^K x_{ik} = p_i, \quad \sum_{k=1}^K q_k = 1, \quad i = 1, \ldots, I
    \end{equation}

    \noindent where $p_i > 0, \sum_{i=1}^I p_i = 1$ represent normalized supply volumes, $x_{ik}$ denote transportation volumes, $d(\xi_i, y_k) = \sqrt{\sum_{j=1}^n | \xi_{ij} - y_{kj} |^2}$ is the metric between elements in the objective function (\ref{transport:eq}), $Y \subset \mathbb{R}^{n}$ is a common constraint set for variables $\{y_k, k = 1, \ldots, K\}$, and $n, I, K \in \mathbb{N}$. We chose the Euclidean norm as a distance metric $d(\xi_i, y_k)$, as the similar colors would have close primary colors magnitude, thus they have small distance in Euclidean space between each other.

    In the context of color quantization, the distribution $\{\xi_i\}$ represents a set of pixels from the original image, while $\{y_k\}$ denotes an optimal color palette with the number of colors $K$ set manually as a hyperparameter. Both $\{\xi_i\}$ and $\{y_k\}$ are subsets of the space $\{(x, y, z) \in \mathbb{Z}^3 : 0 \leq x, y, z \leq 255\}$, representing the combined intensity of three primary colors for each image pixel. Prior to solving the problem (\ref{transport:eq}), we normalize the original pixel distribution $\{\xi_i\}$ to the unit cube $[0, 1]^3 \subset \mathbb{R}^3$.

    \section{Stochastic Color Quantization}

    Paper \cite{Norkin_2024} proposed a solution to the problem (\ref{transport:eq}) by reformulating it as a non-convex non-smooth stochastic transportation problem:

    \begin{eqnarray}
        \label{stochastic-transport:eq}
            &\min_{y = \{ y_1, \ldots, y_K \} \in Y^K \subset \mathbb{R}^{nK}} F(y_1, \ldots, y_k) \nonumber \\
            &F(y) = F(y_1, \ldots, y_k) = \sum_{i=1}^I p_i \min_{1 \leq k \leq K} d(\xi_i, y_k)^r = \nonumber \\
            &= \mathbb{E}_{i \sim p} \min_{1 \leq k \leq K} d(\xi_i, y_k)^r
    \end{eqnarray}
    
    \noindent where $\mathbb{E}_{i \sim p}$ denotes the expected value over a set of image pixels $\{\xi_i\}$ with corresponding probabilities $\{p_i\}$.
    
    To solve the stochastic transportation problem (\ref{stochastic-transport:eq}), we employed the Stochastic Quantization algorithm, as discussed in \cite{Lakshmanan_Pichler_2023} and \cite{Norkin_2024}:
    
    \begin{eqnarray}
        \label{stochastic-quantization:eq}
            &y_k^{(t+1)} = \pi_Y (y_k^{(t)} - \rho g_k^{(t)}), \quad k = 1, \ldots, K, \nonumber \\
            &g_k^{(t)} = \begin{cases}
                r \| \tilde{\xi}^{(t)} - y_k^{(t)} \|^{r - 2} (y_k^{(t)} - \tilde{\xi}^{(t)}), & k = k^{(t)}, \\
                0, & k \neq k^{(t)}
            \end{cases} \nonumber \\
            &k^{(t)} \in S(\tilde{\xi}^{(t)}, y^{(t)}) = \argmin_{1 \leq k \leq K} d(\tilde{\xi}^{(t)}, y_k^{(t)}), \quad t = 0, 1, \ldots
    \end{eqnarray}

    \noindent where $y^{(t)} = (y_1^{(t)}, \ldots, y_K^{(t)})$, $\tilde{\xi}^{(t)}$ is a random element from the set $\{ \xi_i, i = 1, \ldots, I \}$, $\rho_t > 0$ is a learning rate parameter, $\pi_Y$ is the projection operator onto the set $Y \subseteq \mathbb{R}^n$, $r\geq3$ is the degree of the norm, $K$ is the number of colors in the optimal palette, and $t = 0, 1, \ldots$ is the number of iterations. To calculate the update for the element $y_k^{(t)}$ based on the gradient value $g_k^{(t)}$ we used Stochastic Gradient Descent (SGD) algorithm \cite{ermoliev1976stochastic} and its non-convex non-smooth extension from \cite{Ermolev_Norkin_1998}.

    In contrast to traditional clustering algorithms such as K-Means \cite{Lloyd_1982}, which update all components $y_k^{(t)}$ and uses all elements $\{\xi_i\}$ per iteration $t$, our approach processes only a single element of $y^{(t)}$ and uses one random element $\tilde{\xi} \in \{ \xi_1, \ldots, \xi_I \}$ per iteration. This eliminates memory constraints and enhances scalability for high-resolution images. Furthermore, the paper \cite{Norkin_2024} established local convergence conditions for the algorithm (\ref{stochastic-quantization:eq}) under specific initial conditions, which is not typically guaranteed for traditional algorithms.

    Given the non-convex nature of the problem, the initialization (seeding) of elements $y^{(0)} = \{ y_k, k = 1, \ldots, K \}$ plays a crucial role in the algorithm's convergence. While uniform sampling from the original distribution $\Xi$ is the fastest seeding approach, it can compromise accuracy. We implemented the seeding technique proposed in \cite{Arthur_Vassilvitskii_2007}, where the selection probability of a new element in $y_k^{(0)}$ is proportional to its distance from previously selected elements $\{ y_k^{(0)} \}$. The seeding algorithm proceeds as follows:
    
    \begin{enumerate}
        \item Sample the initial element $y_0^{(0)}$ uniformly from $\Xi$
        \item Select the next element $y_k^{(0)}$ with probability $\frac{\min_{1 \leq s \leq k} \| \xi_j - y_s^{(0)} \|^2}{\sum_{i=1}^I \min_{1 \leq s \leq k} \| \xi_i - y_s^{(0)} \|^2}$
        \item Repeat step 2 until $y^{(0)}$ contains $K$ elements
    \end{enumerate}
    
    The research \cite{Norkin_2024} empirically demonstrates the efficiency of this seeding technique and compares the convergence rates of modified versions of the algorithm (\ref{stochastic-quantization:eq}) with adaptive learning rates.

    \section{Numerical Experiments}

    To evaluate the efficiency of Stochastic Color Quantization, we conducted lossy compression experiments on a subset of images from the ImageNet dataset \cite{Russakovsky_2015}, varying in resolution and stored in JPEG format. The algorithm was implemented in Python, utilizing NumPy \cite{harris2020array} for efficient CPU-based tensor operations. For loading, processing and storing images from the file system we utilized Pillow library \cite{Andrew_Murray_2024}. All experiments were executed on a virtual machine with a dual CPU core and 8GB RAM. For reproducibility, we set the random seed to 42 and used consistent hyperparameter values ($\rho=0.001$, $r=3$) across all experiments. The source code and experimental results are publicly available in our GitHub repository \cite{Kozyriev_2024}.

\begin{figure}
        \centering
        \includegraphics[width=\textwidth]{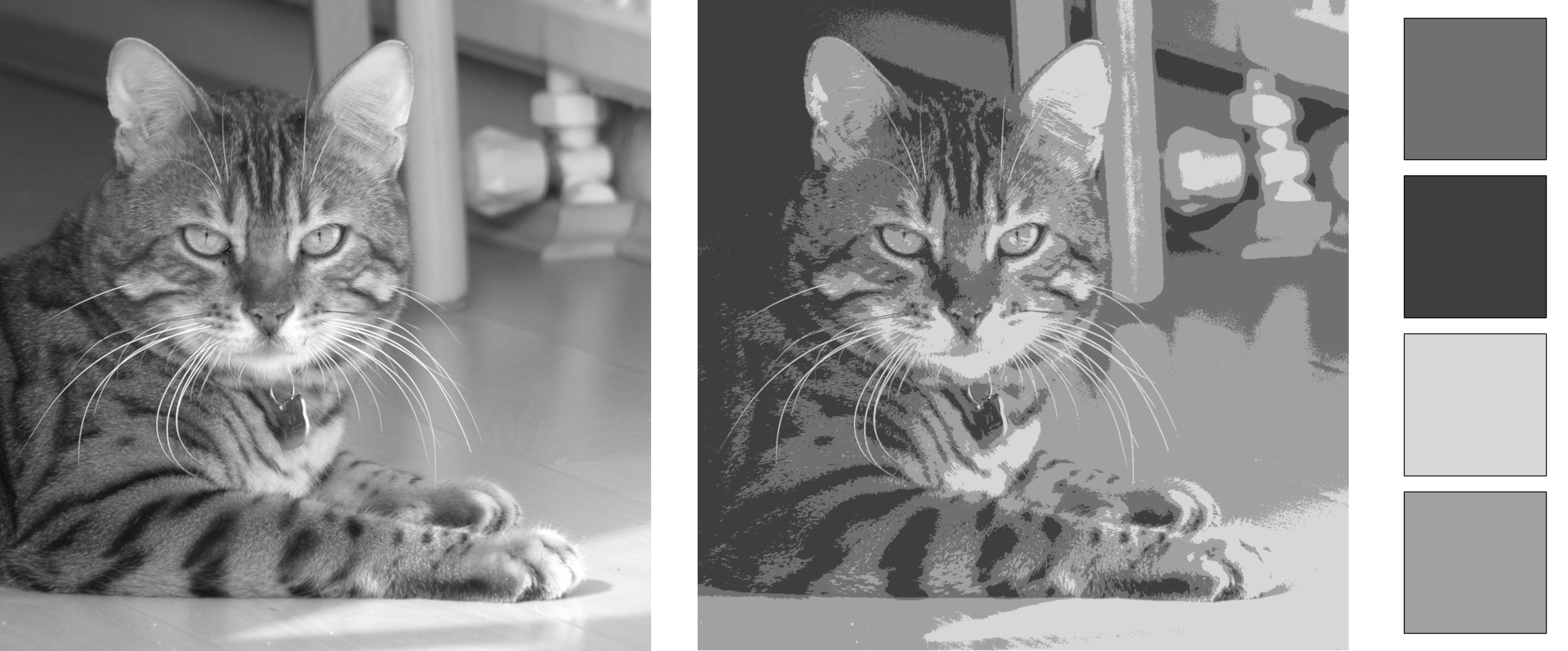}
        \caption{Original image (left) compared with compressed image with lower quality and an optimal color palette with 4 colors (right).}
        \label{monochrome:fig}
    \end{figure}

    We utilized the Mean Squared Error (MSE) as our primary distortion metric to quantify the similarity between original and compressed images across various resolutions and optimal color palette sizes (see Table \ref{mse:table}). To visually demonstrate the algorithm’s performance, we applied it to a monochrome image labeled \textit{ILSVRC2012\_val\_00023267} from the dataset. The original image (1200$\times$1206 resolution, 248 distinct colors, 790.0KB) was compressed to 207.4KB with an optimal palette of four colors: \texttt{\#707070}, \texttt{\#3e3e3e}, \texttt{\#d7d7d7}, and \texttt{\#a1a1a1} (see Fig.~\ref{monochrome:fig}). The algorithm converged with a transport value of $F(y) = 129465.2$ and a distortion of $MSE = 0.0037$.

    \begin{table}
        \caption{Comparison of the distortion metric (MSE) between original and compressed images for different resolutions.}
        \label{mse:table}
        \begin{tabularx}{\textwidth}{|c|*{3}{Y|}} % Use 'Y' instead of 'X'
            \hline
            \multirow{2}{*}{Optimal color palette size (K)} & \multicolumn{3}{c|}{Image resolution} \\
            \cline{2-4}
            & 438x500 & 1200x1206 & 1606x2400 \\
            \hline
            4  & $MSE = 0.0043$ & $MSE = 0.0056$ & $MSE = 0.0049$ \\
            8  & $MSE = 0.0019$ & $MSE = 0.0029$ & $MSE = 0.0019$ \\
            12 & $MSE = 0.0013$ & $MSE = 0.0018$ & $MSE = 0.0012$ \\
            24 & $MSE = 0.0009$ & $MSE = 0.0010$ & $MSE = 0.0005$ \\
            36 & $MSE = 0.0007$ & $MSE = 0.0007$ & $MSE = 0.0004$ \\
            \hline
        \end{tabularx}
    \end{table}

    \section{Conclusions}

    This study introduces a scalable algorithm for solving the color quantization problem without memory constraints, demonstrating its efficiency on a subset of images from the ImageNet dataset \cite{Russakovsky_2015}. The convergence speed of the algorithm can be further enhanced by modifying the update rule (\ref{stochastic-quantization:eq}) with alternative methods to Stochastic Gradient Descent (SGD) that incorporate adaptive learning rates, as explored in \cite{Norkin_2024}. Moreover, the stochastic nature of the proposed solution enables the utilization of parallelization techniques to simultaneously update the positions of multiple quants $\{y_k, k = 1, \ldots, K\}$, potentially leading to significant performance improvements. This aspect of parallelization and its impact on algorithm efficiency presents a topic for future research. The proposed method not only addresses the limitations of existing color quantization techniques but also opens up new possibilities for optimizing image compression algorithms in resource-constrained environments.

    \bibliographystyle{splncs04}
    \bibliography{references}
\end{document}